\documentclass{amsart}

\usepackage{amssymb}
\usepackage{amsfonts}
\usepackage{amsmath}
\usepackage{mathrsfs}
\newtheorem{theorem}{Theorem}
\theoremstyle{plain}

\newtheorem{corollary}{Corollary}

\newtheorem{definition}{Definition}

\newtheorem{lemma}{Lemma}

\newtheorem{proposition}{Proposition}
\newtheorem{remark}{Remark}

\numberwithin{equation}{section}
\newcommand{\Prob}{\mathbb{P}} 
\newcommand{\Pba}{\mathbb{P}_{ba}} 
\newcommand{\B}{\mathfrak{B}} 
\newcommand{\A}{\mathscr{A}} 

\newcommand{\R}{\mathbb{R}} 
\newcommand{\N}{\mathbb{N}} 

\newcommand{\restr}[2]{#1\vert #2}

\newcommand{\abs}[1]{\vert #1\vert} 
 
\newcommand{\norm}[1]{\Vert #1\Vert} 
\newcommand{\dnorm}[1]{\left\Vert #1\right\Vert}

\newcommand{\net}[3]{\left\langle #1_{#2}\right\rangle_{#2\in #3} } 
\newcommand{\nnet}[3]{\left\langle #1\right\rangle_{#2\in #3} } 
\newcommand{\seq}[2]{\net{#1}{#2}{\mathbb{N}}} 
\newcommand{\sseq}[2]{\nnet{#1}{#2}{\mathbb{N}}} 
\newcommand{\seqn}[1]{\seq{#1}{n}}

\newcommand{\K}{\mathcal{K}}
\newcommand{\M}{\mathscr{M}}
\newcommand{\La}{\mathscr{L}}
\newcommand{\T}{\mathcal{T}} 

\newcommand{\Gg}{\Gamma(\gamma)}
\newcommand{\set}[1]{\mathbf{1}_{#1}}

\begin{document}
\title[Sure wins]{Sure Wins, Separating Probabilities and the 
Representation of Linear Functionals}
\author{Gianluca Cassese}
\address{Universit\`{a} Milano Bicocca and University of Lugano}
\email{gianluca.cassese@lu.unisi.ch}
\curraddr{Department of Statistics, Building U7, Room 239, via Bicocca 
degli Arcimboldi 8, 20126 Milano - Italy}
\date
\today


\begin{abstract}
We discuss conditions under which a convex cone $\K\subset\R^\Omega$ 
admits a finitely additive probability $m$ such that $\sup_{k\in\K}m(k)\leq0$. 
Based on these, we characterize those linear functionals that are representable
as finitely additive expectations. A version of Riesz decomposition 
based on this property is obtained as well as a characterisation of 
positive functionals on the space of integrable functions. 
\end{abstract}

\maketitle

\section{Introduction}

A long standing approach to probability, originating from the
seminal work of de Finetti, views set functions $P$ as maps which
assign to each set (event) $E$ in some class $\A$ the price $P(E)$
for betting 1 dollar on the occurrence of $E$. A set function
generating a betting system which admits no sure wins was termed
coherent by de Finetti who proved in \cite{definetti} that a set
function on a finite algebra $\A$ is coherent if and only if it is
a probability. Since then this result has been extended and generalized 
by various authors, among which Heath and Sudderth \cite{heath
sudderth}, Lane and Sudderth \cite{lane sudderth} and Regazzini
\cite{regazzini}, to name but a few; Borkar et al. \cite{borkar}
is a more recent example. 

In this paper we examine the absence of sure wins for a convex cone 
$\K$ of real valued functions on some arbitrary set $\Omega$, obtaining 
conditions for the existence of a finitely additive probability measure 
$m$ such that $\sup_{k\in\K}m(k)\leq0$, i.e. a \textit{separating probability}. 
The special case in which $\K$ is the kernel of some linear 
functional leads to the characterization of those functionals
that admit the representation as finitely additive expectations, a topic 
addressed by Berti and Rigo in a highly influential paper \cite{berti rigo}. 
A version of Riesz decomposition based on this representation property is 
obtained. 

Throughout the paper $\Omega$ will be a fixed set, $2^{\Omega}$ its 
power set, $\R^\Omega$ and $\B$ the classes of real-valued and bounded 
functions on $\Omega$ respectively (the latter endowed with the topology 
induced by the supremum norm). All spaces of real-valued functions on 
$\Omega$ (e.g. bounded or integrable) will be considered as equipped with 
pointwise ordering, with no further mention. The lattice notation $f^+$ and 
$f^-$ wll be used to denote the positive and negative parts of $f\in\R^\Omega$. 
The term probability is used to designate positive, finitely additive set 
functions $m$ on $2^{\Omega}$ (in symbols, $m\in ba_+$) such that $m(\Omega)=1$. 
The symbol $\Pba$ will be used to denote the family of all probability measures; 
$\Prob$ the subfamily of all countably additive probability measures. If 
$\A\subset 2^\Omega$ then by $\mathscr{S}(\A)$ and $\B(\A)$ we denote the 
class of simple functions generated by $\A$ and its closure in $\B$. We 
adopt the useful convention of identifying single-valued functions with their range 
so that, for example, we may use $1$ either to denote an element of $\R$, or a function 
$f$ on $\Omega$ such that $f(\omega)=1$ for all $\omega\in\Omega$. In 
the terminology adopted throughout the following sections a \textit{sure win} 
is defined to be an element of $\R^\Omega$ which exceeds $1$.

We recall that $f\in\R^\Omega_+$ is integrable with respect to $m\in ba_+$, 
in symbols $f\in L(m)$, if and only if 
\begin{equation}
\sup\left\{m(h):h\in\B, 0\leq h\leq f\right\}<\infty
\label{expectation}
\end{equation}
The integral $m(f)$ coincides then with the left hand side of (\ref{expectation}); 
moreover, $f\wedge n$ converges to $f$ in $L(m)$ \cite[theorem III.3.6]{bible}.
A special notion of convergence in $L(m)$ will be used in the following. A 
sequence $\seqn{f}$ is said to converge orderly in $L(m)$ to $f$ if $f_n\in L(m)$ 
for all $n$ and there exists a pointwise decreasing sequence $\seqn{\bar f}$ in 
$L(m)_+$ which converges to $0$ in $L(m)$ and is such that $\abs{f_n-f}\leq\bar f_n$ 
for $n\geq 1$. It is easily seen that if a sequence $\seqn{f}$ converges to $f$
orderly in $L(m)$ then so does any of each subsequences; moreover, the space
of sequences converging orderly in $L(m)$ is a vector space.

\section{Separating Probabilities}

Fix a convex cone $\K\subset\R^\Omega$ (that is 
$f+g,\ \lambda f\in\K$ whenever $f,g\in\K$ and $\lambda\geq 0$) 
and let $\K_b=\left\lbrace k\in\K:k^-\in\B\right\rbrace$. For each $f\in\R^\Omega$ 
let $\mathcal{U}(f)=\{\alpha\in\R:\alpha+k\geq f\text{ for some }k\in\K\}$ 
and define $\pi_\K :\R^\Omega\rightarrow\overline{\R}$ as
\begin{equation}
\pi_\K(f) =\inf\{\alpha :\alpha\in\mathcal{U}(f)\} %
\footnote{The functional $\pi_\K $ is well known 
in mathematical finance under the name of \textit{superhedging price}.}
\label{hedge}
\end{equation}
From (\ref{hedge}), $\pi_\K $ is monotonic, $\pi_\K(\lambda +f) =\lambda +\pi_\K(f)$ 
for each $\lambda\in\R$ and $f\in\R^\Omega$ and $\pi_\K(f)\leq\sup_{\omega\in\Omega}f(\omega)$ 
(as $0\in\K$). Since $\K$ is a convex cone, $\mathcal{U}(f)+\mathcal{U}(g)\subset\mathcal{U}(f+g)$
and $\mathcal{U}(\lambda f)=\lambda\mathcal{U}(f)$ for $\lambda>0$: $\pi_\K$ is 
thus subadditive and positively homogeneous; moreover, $\pi_\K(k)\leq 0$ for all 
$k\in\K$.

Given that $\pi_\K(0)=2\pi_\K(0)\leq 0$ and $\pi_\K(1) =\pi_\K(0)+1$, then $\pi_\K(0) >-\infty$ 
implies $\pi_\K(0)=0$ and $\pi_\K(1)=1$. Moreover there is $k\in\K$ such that $k\geq 1$ 
if and only if $\pi_\K(1)\leq0$. Thus:

\begin{lemma}
\label{lemma superhedge}
Let $\K\subset\R^\Omega$ be a convex cone. Then the following are 
equivalent: (i) $\pi_\K(0)> -\infty $, (ii) $\pi_\K(0) =0$, (iii) $\pi_\K (1) =1$, (iv) $\K$ 
contains no sure wins.
\end{lemma}

Denote $L(\pi_\K) =\left\{ f\in\R^\Omega:\pi_\K(\abs{f})<\infty\right\} $. 
It is clear that $\B\subset L(\pi_\K)$. Define also
\begin{equation}
\label{M(K)}
\M(\K)=\left\{ m\in\Pba: \K\subset L(m),\quad\sup_{k\in\K}m(k)\leq
0\right\}
\end{equation}
and let $\M(\K_b)$ be defined likewise. We shall refer to elements of 
$\M(\K)$ as \textit{separating probabilities} for $\K$. It is clear that 
if $m\in\M(\K_b)$ then $L(\pi_\K)\subset L(m)$.

\begin{proposition}
\label{proposition M(Kb)}
Let $\K\subset\R ^{\Omega}$ be a convex cone. Then $\M(\K_b)$ is non empty if 
and only if $\K$ contains no sure wins.
\end{proposition}

\begin{proof}
Assume that $\K$ contains no sure wins. By Lemma \ref{lemma superhedge} 
and the Hahn Banach Theorem, we may find a linear functional $\phi $ on 
$\B$ such that $\phi\leq\pi_\K$ on $\B$ and $\phi(1) =1$. If $f\in\B_+$ then 
$\phi(f) =-\phi(-f)\geq-\pi_\K(-f)\geq 0$. Therefore $\phi$ is positive and, 
since continuous \cite[V.2.7]{bible}, it may be represented as the expectation 
with respect to some $m\in\Pba$. If $f\in L(\pi_\K)_+$, the left hand side of 
(\ref{expectation}) is bounded by $\pi_\K(f)$ so that $L(\pi_\K)\subset L(m)$. 
Then $\K_b\subset L(m)$ and $$m(k)=\lim_nm(k\wedge n)\leq\pi_\K(k)\leq 0\qquad k\in\K_b$$
so that $m\in\M(\K_b)$. If $m\in\M(\K_b)$ and $k\in\K$ is a sure win, then 
$k\in\K_{b}$ and $m(k)\leq 0$, a contradiction.
\end{proof}

A classical application of Proposition \ref{proposition M(Kb)} considers the 
collection $\K$ of all finte sums of the form $\sum_na_n(\set{F_n}-\lambda(F_n))$ 
where $a_1,\ldots,a_N$ are real numbers, $F_1,\ldots,F_N$ are elements of some 
$\A\subset 2^{\Omega}$ and $\lambda:\A\rightarrow\R$. It is then clear that 
$\K$ admits no sure wins if and only if there is $m\in\Pba$ such that 
$\restr{m}{\A}=\lambda$. If the sums in $\K$ are allowed to admit countably 
many terms provided $\sum_n\abs{a_n\lambda(F_n)}<\infty$, then $m$ will possess 
the additional property that $m(\bigcup_nF_n)=\sum_nm(F_n)$ when $\seqn{F}$ is 
a disjoint sequence in $\A$. This informal statement is essentially a 
reformulation of \cite[theorems 5 and 6, p. 2074]{heath sudderth}\footnote{However  
we do not restrict $\A$ nor $\lambda$. Heath and Sudderth seem to suggest that 
the existence of $m$ need not exclude sure wins while it is clear that this cannot 
be the case. A less general version of this result was also proved, with different 
methods, in \cite[theorem 2, p. 420]{borkar}}. It admits an interesting generalisation
to the case of concave integrals, a special case of the monotone integral of Choquet 
treated, e.g., in \cite{greco}. 

\begin{definition}
\label{definition concave integral}
An extended real-valued functional $\gamma$ on a convex cone $\La\subset\R^\Omega$ 
is a concave integral if it is positively homogeneous, monotone, superadditive and 
such that $\gamma(c+f)=\gamma(c)+\gamma(f)$ when $c,f\in\La$ and $c$ is a constant. 
\end{definition} 

If $\gamma$ is a concave integral on $\La$ we define its core to be the set
\begin{equation}
\label{core}
\Gg=\{\lambda\in\ ba_+:\ \La\subset L(\lambda),\ \gamma(f)\leq\lambda(f),\ f\in\La\}
\end{equation} 

The following Lemma is essentially a restatement of a result of Shapley \cite[theorem 2, 
p. 18]{shapley}. It characterises the properties of a concave integral in terms of its core.

\begin{lemma}
\label{lemma shapley}
Let $\La\subset\R^\Omega$ be a convex cone that contains the constants and is such that
$f\in\La$ implies $f^+\in\B$. Let $\gamma:\La\rightarrow\R$ be a concave integral and
$\gamma(1)>0$. Then $\gamma(1)<\infty$ if and only if for each convex set 
$C\subset\La\cap\B$ such that $\gamma(C)\equiv\sup_{f\in C}\gamma(f)<\infty$ there 
exists $\lambda_C\in\Gg$ such that 
\begin{equation}
\label{infsup}
\sup_{f\in C}\lambda_C(f)=\gamma(C)
\end{equation}
\end{lemma}

\begin{proof}
Assume, upon normalization, $\gamma(1)=1$ and suppose that
\begin{equation}
\alpha(k-\gamma(C))\geq 1+\sum_{n=1}^N(f_n-\gamma(f_n))
\label{sure win}
\end{equation}
for some choice of $\alpha\geq 0$, $k\in C$ and $f_n\in\La$, $n=1,\ldots,N$. 
The value under $\gamma$ of the left hand side of (\ref{sure win}) is less 
than $0$ while that of the right hand side exceeds $1$, contradicting monotonicity. 
Thus the collection $\K_C$ of finite sums of the form 
$
\sum_{1\leq n\leq N}(\gamma(f_n)-f_n)+\alpha(k-\gamma(C))
$
for $\alpha$, $k$ and $f_n$, $n=1,\ldots,N$ as above contains no sure win; 
moreover, it is a convex cone of uniformly lower bounded functions on $\Omega$. 
According to Proposition \ref{proposition M(Kb)}, there exists 
$\lambda_C\in\M(\K_C)$: thus, $\lambda_C(f)\geq\gamma(f)$ for each $f\in\La$ 
(i.e. $\lambda_C\in\Gamma(\gamma)$) and $\lambda_C(k)\leq\gamma(C)$ whenever 
$k\in C$, proving (\ref{infsup}). The converse is obvious. 
\end{proof}
Lemma \ref{lemma shapley} has an interesting implication.

\begin{corollary}
\label{corollary common extension}
Let $\T$ be a collection of subsets of some set $T$, with $\{T\}=\tau_0\in\T$. 
For each $\tau\in\T$, let $\La_\tau$ be a linear subspace of $\B$ with $1\in\La_{\tau_0}$ 
and $\phi_\tau$ a linear functional on $\La_\tau$. The following are equivalent:
\begin{enumerate}
\item[(\textit{i})] the collection $(\phi_\tau:\tau\in\T)$
is coherent in the sense that\footnote{The inequality that follows is meant to hold 
pointwise in $\Omega\times T$} 
$$\sup\left\{\sum_{n=1}^N\phi_{\tau_n}(b_n):\ 
                       b_n\in\La_{\tau_n},\ 
                       \sum_{n=1}^Nb_n\set{\tau_n}\leq 1,\ N\in\N\right\}<\infty$$
\item[(\textit{ii})] there exists $\lambda\in ba(\Omega\times T)$ such that
$\lambda(b\set{\tau})=\phi_\tau(b)$ for each $b\in\La_\tau$ and $\tau\in\T$
\end{enumerate}
\end{corollary}
\begin{proof}
Assume (\textit{i}) and define the functional $\gamma$ on $\B(\Omega\times T)$ implicitly as
\begin{equation}
\label{generated integral}
\gamma(b)=\sup\left\{\sum_{n=1}^N\phi_{\tau_n}(b_n):\ 
                       b_n\in\La_{\tau_n},\ 
                       \sum_{n=1}^Nb_n\set{\tau_n}\leq b,\ N\in\N\right\}
\end{equation}
It is readily seen that $\gamma$ is monotone, superadditive and positively 
homogeneous. (\textit{i}) implies that $\gamma(1)<\infty$ and that $\gamma$
is real-valued while $1\in\La_{\tau_0}$ implies that $\gamma$ is additive 
relative to the constants. (\textit{ii}) 
follows from (\textit{i}), Lemma \ref{lemma shapley} and the fact that each 
$\La_\tau$ is a linear space: simply choose $\lambda\in\Gg$. If $\lambda$ is 
as in (\textit{ii}) and $\sum_{n=1}^Nb_n\set{\tau_n}\leq 1$ where $b_n\in\La_{\tau_n}$ 
$n=1,\ldots,N$ then $\sum_{n=1}^N\phi_{\tau_n}(b_n)=\lambda\left(\sum_{n=1}^Nb_n%
\set{\tau_n}\right)\leq\norm{\lambda}$.
\end{proof}

\begin{remark} 
Writing $\tau\leq\upsilon$ when $\tau\subset\upsilon$ makes of course $\T$ 
into a partially ordered set. If $(\phi_\tau:\tau\in\T)$ is coherent in the 
sense of Corollary \ref{corollary common extension} and if $(\La_\tau:\tau\in\T)$ 
is increasing in $\tau$ then necessarily $\restr{\phi_\upsilon}{\La_\tau}\geq\phi_\tau$
whenever $\tau,\upsilon\in\T$ and $\tau\leq\upsilon$. This conclusion has a 
direct application to the theory of finitely additive supermartingales, treated 
in \cite{general}.
\end{remark}

Much of this section rests on the conclusion, established in Proposition 
\ref{proposition M(Kb)}, that $\K_b$ admits a separating probability in the 
absence of sure wins. This result, however, does not have an extension to $\K$ of a 
corresponding simplicity. To this end we shall need some results on the 
representation of linear functionals, to be developed in the next section. 

\section{The Representation of Linear Functionals}

It is the purpose of this section to establish conditions for a linear 
functional $\phi$ on some linear subspace $\La$ of $\R^\Omega$ with 
$1\in\La$ to admit the representation 
\begin{equation}
\label{representation}
\phi(f)=\phi(1)m(f)\qquad f\in\La
\end{equation}
for some $m\in ba$ such that $\La\subset L(m)$, referred to as a 
\textit{representing measure} for $\phi$. We use the symbols $\K^{\phi}$ 
and $\K^{\phi}_b$ to denote the sets $\{f\in\La:\phi(f)=0\}$ and 
$\{f\in\K^{\phi}:f^-\in\B\}$, respectively. If $\phi(1)\neq0$, then
 $\K^{\phi}_b=\{f-\phi(1)^{-1}\phi(f):f\in\La,f^-\in\B\}$. 
Thus if $\La$ is a vector sublattice of $\R^\Omega$ then $m\in\M(\K^{\phi}_b)$ 
implies $\La\subset L(m)$ and $\phi(f)=\phi(1)m(f)$ 
for every $f\in\La\cap\B$ (which clarifies the connection between 
separating probabilities and representing measures).

The content of this section, as will soon become clear, owes much to the work of 
Berti and Rigo \cite{berti rigo}. 

\begin{theorem}
\label{theorem extension} 
Let $\A\subset 2^{\Omega}$ be an algebra, $\mu\in ba(\A)$, $\La$ 
a vector sublattice of $L(\mu)$ with $1\in\La$ and $\phi$ a positive linear 
functional on $\La$. Denote by $\La^*$ the set of limit points of sequences
from $\La$ which converge orderly in $L(\mu)$. The following are equivalent.
\begin{enumerate}
\item[(\textit{i})] $\phi$ extends to a monotone function $\phi^*:\La^*\rightarrow\R$;
\item[(\textit{ii})] $\lim_n\phi(h_n)=0$ whenever $\seqn{h}$ is a sequence in
$\La$ which converges to $0$ orderly in $L(\mu)$;
\item[(\textit{iii})] $-\infty<\lim_n\phi(g_n)\leq\lim_n\phi(f_n)<\infty$ whenever $\seqn{f}$ 
and $\seqn{g}$ are sequences in $\La$ which converge orderly in $L(\mu)$ to $f$ 
and $g$ respectively, with $f\geq g$;
\item[(\textit{iv})] $\phi$ admits a positive representing measure $m$ such that 
$m^*(h)\equiv\lim_nm(h_n)$ exists in $\R$ and is unique for every sequence 
$\seqn{h}$ in $\La$ which converges to $h$ orderly in $L(\mu)$.
\end{enumerate}
Moreover, if $\phi$ is a positive linear functional on a vector sublattice $\La$ of 
$\R^{\Omega}$ with $1\in\La$ then there exists a unique positive linear functional 
$\phi^\perp$ on $\La$ such that $\phi^{\perp}(1)=0$ and that
\begin{equation}
\phi(f)=\phi(1)m(f)+\phi^\perp(f)\qquad f\in\La
\label{riesz}
\end{equation}
for some $m\in ba_+$ satisfying $\La\subset L(m)$. 
\end{theorem} 

\begin{proof}
Let $\seqn{h}$ be as in (\textit{ii}) and $\seqn{\bar h}$ be 
a decreasing sequence in $L(m)$ converging to $0$ in $L(m)$ and 
such that $\bar h_n\geq\abs{h_n}$, $n=1,2,\ldots$. Fix a sequence $\seqn{\alpha}$ 
in $\R_+$ such that $\lim_n\alpha_n=\infty$. Any subsequence of $\seqn{h}$ 
admits a further subsequence (still denoted by $\seqn{h}$ for convenience) 
such that $\sum_n\alpha_n\norm{h_n}<\infty$. Fix $\eta>0$ arbitrarily and set 
\begin{equation}
\label{h^eta}
h_n^\eta=(h_n-\eta)^+,\ g^\eta_k=\sum_{n\leq k}\alpha_nh_n^\eta\quad\text{and}
\quad g^\eta=\sum_n\alpha_nh_n^\eta
\end{equation}
Then, $\left\{\sum_{n>k}\alpha_nh_n^\eta>\epsilon\right\}\subset\{\bar h_k\geq\eta\}$ and 
$\dnorm{\sum_{k<n\leq k+p}\alpha_nh_n^\eta}\leq\sum_{n>k}\alpha_n\norm{h_n^\eta}%
\leq\sum_{n>k}\alpha_n\norm{h_n}$. Thus, $\seq{g^\eta}{k}$ is an increasing sequence 
in $\La$ which converges orderly in $L(\mu)$ to $g^\eta\in\La^*$ \cite[theorem III.3.6]{bible}. If 
(\textit{i}) holds then 
$\alpha_n^{-1}\phi^*(g^\eta)\geq\phi(h_n^\eta)\geq\phi(h_n)-\eta\phi(1)$ 
so that $\lim_n\phi(h_n)=0$, i.e. (\textit{ii}) holds as well. 
Let $\seqn{g}$ and $\seqn{f}$ be as in (\textit{iii}). The inequality 
$f_n-g_n\geq (f_n-f)+(g-g_n)$ together with (\textit{ii}) induces the 
conclusion that $(f_n-g_n)^-$ converges to $0$ orderly in $L(\mu)$ and thus that 
$\liminf_n\phi(f_n)=\liminf_n\{\phi(g_n)+\phi((f_n-g_n)^+)\}\geq\liminf_n\phi(g_n)$. 
The case in which $\seqn{g}$ is a subsequence of $\seqn{f}$ suggests that 
$\liminf_n\phi(f_n)=\limsup_n\phi(f_n)$. If $\lim_n\phi(f_n)=\infty$ then 
one may select a subsequence $\sseq{f_{n_k}}{k}$ such that, letting 
$h_k=f_{n_{k+1}}-f_{n_k}$, $\lim_k\phi(h_k)=\infty$. However this contrasts with 
(\textit{ii}) since the sequence $\seq{h}{k}$ converges to $0$ orderly in
$L(\mu)$. This proves (\textit{iii}). 
In the general case in which $\La$ is a vector sublattice of $\R^\Omega$, fix 
$f\in\La_+$ and choose $m\in\M(\K^\phi_b)$ if $\phi(1)>0$, or $m=0$ otherwise. Then, 
\begin{equation}
\label{density}
\phi(f)=\lim_n\phi(f\wedge n)+\lim_n\phi(f-(f\wedge n))=\phi(1)m(f)+\phi^{\perp}(f)
\end{equation}
a conclusion which extends to general $f\in\La$ by considering 
$f^+$ and $f^-$ separately. The functional $\phi^{\perp}$, as 
defined implicitly in (\ref{density}), is clearly positive, linear 
and such that $\phi^{\perp}(1)=0$. Decomposition (\ref{riesz}) 
thus exists. If $\phi(f)=\phi(1)v(f)+\psi^{\perp}(f)$ were another 
decomposition such as (\ref{riesz}), with $v\in ba_+$, $\La\subset L(v)$ 
and $\psi^{\perp}$ a positive, linear functional on $\La$ with 
$\psi^{\perp}(1)=0$, then $f\in\La_+$ would imply
$$
(\phi^{\perp}-\psi^{\perp})(f)=\lim_n(\phi^{\perp}-\psi^{\perp})(f-(f\wedge n))
=\phi(1)\lim_n(m+v)(f-(f\wedge n))=0
$$
which proves uniqueness of (\ref{riesz}). Returning to the case $\La\subset L(\mu)$,
if (\textit{iii}) holds, then it is obvious from (\ref{density}) that $\phi^{\perp}=0$; 
in addition the limit $\lim_nm(h_n)$ exists in $\R$ for each sequence $\seqn{h}$ in
$\La$ which converges orderly in $L(\mu)$ and does not depend but on the limit point $h$. 
\end{proof}

One noteworthy implication of Theorem \ref{theorem extension}, obtained 
by replacing $\La$ with $L(\mu)$, is the following

\begin{theorem}
\label{theorem L(mu)}
Let $\A\subset 2^{\Omega}$ be an algebra and $\mu\in ba(\A)$. Every 
positive linear functional $\phi$ on $L(\mu)$ admits a positive representing 
measure $m$ such that $\lim_nm(h_n)=0$ for every sequence $\seqn{h}$ in $L(\mu)$
which converges to $0$ orderly in $L(\mu)$. 
\end{theorem}

Given that $L(\mu)$ is a normed Riesz space, its dual space is a 
vector lattice \cite[theorem 12.1, p. 175]{aliprantis}. Thus Theorem 
\ref{theorem L(mu)} also implies that continuous linear functionals, 
decomposing as the difference of two positive linear functionals, 
admit a representing measure \cite[theorem 7, p. 3255]{berti rigo}. 

Another application concerns more general functionals. In fact it is 
clear that the implication (\textit{i})$\rightarrow$(\textit{ii}) in 
Theorem \ref{theorem L(mu)} does not require $\phi$ to be linear. 

\begin{theorem}
\label{theorem Choquet}
Let $\La\subset\R^\Omega$ be either (\textit{i}) a Banach lattice containing 
the constants or (\textit{ii}) $\La=L(\mu)$ for some $\mu\in ba(\A)$ and some 
algebra $\A\subset2^\Omega$. Assume that $\phi:\La\rightarrow\R$ is a monotone 
functional such that 
\begin{equation}
\label{pseudo homogeneous}
\lim_n\inf_{\{f\in\La:\phi(f)>\eta\}}\phi(nf)=\infty\qquad\eta>0
\end{equation}
and, under (\textit{ii}), 
\begin{equation}
\label{pseudo additive}
\lim_{k\downarrow 0}\sup_{f\in\La}\{\phi(f)-\phi(f-k)\}=0
\end{equation}
Then, $\limsup_n\phi(h_n)\leq 0$ when $\seqn{h}$ converges 
to $0$ in norm or, under (\textit{ii}), orderly in $L(\mu)$. In
particular, convex, monotone functionals on a Banach lattice are 
continuous.
\end{theorem}
\begin{proof}
Each subsequence of $\seqn{h}$ contains a further subsequence for which 
it is possible to define $g^\eta_k$ and $g^\eta$ as in (\ref{h^eta}). 
Under (\textit{i}), $\sseq{g^\eta}{k}$ converges to $g^\eta$ in norm
for all $\eta\geq 0$; under (\textit{ii}) only for $\eta>0$. In either case
we conclude that $\phi(g^\eta)\geq\phi(\alpha_nh^\eta_n)\geq\phi(\alpha_n(h_n-\eta))$ 
and, given (\ref{pseudo homogeneous}), $\liminf_n\phi(h_n-\eta)\leq 0$.
Choosing $\eta=0$ under (\textit{i}) or exploiting (\ref{pseudo additive})
under (\textit{ii}) and recalling that the intial choice of the subsequence 
was arbitrary, we conclude that $\limsup_n\phi(h_n)\leq 0$. It is clear that
a convex functional $\phi$ meets (\ref{pseudo homogeneous}), (\ref{pseudo additive})
and, by monotonicity, $\abs{\phi(h)-\phi(h_n)}\leq\phi(\abs{h_n-h})$.
\end{proof}

Given the preceding results, it is now easy to extend Proposition 
\ref{proposition M(Kb)} to $\K$.

\begin{corollary}
\label{corollary M(K)} 
Let $\K\subset\R^\Omega$ be a convex cone. 
Then $\M(\K)$ is non empty if and only if there exist an algebra 
$\A\subset 2^{\Omega}$ and $\mu\in\Pba(\A)$ such that $\K\subset L(\mu)$ 
and that the closure $\overline{C}^{\mu}$ of $C=\K-\mathscr{S}(\A)_+$ 
in the norm topology of $L(\mu)$ admits no sure wins. 
\end{corollary}

\begin{proof}
If $\mu\in\M(\K)$ then $\mu$ is a separating measure for $\overline{C}^\mu$ 
which rules out sure wins. As for sufficiency, observe that ordinary separation 
theorems imply the existence of a continuous linear functional 
$\phi:L(\mu)\rightarrow\R$ such that $\sup_{f\in\overline{C}^\mu}\phi(f)\leq 0$ 
and $1=\phi(1)$. Given that $\K$ contains the origin, $-\mathscr{S}(\A)_+\subset C$ 
so that $\phi$ is positive on $\mathscr{S}(\A)$ and, since $\mathscr{S}(\A)_+$ 
is dense in $L(\mu)_+$ and $\phi$ is $L(\mu)$ continuous, it is positive over 
the whole of $L(\mu)$. The claim follows from Theorem \ref{theorem L(mu)}.
\end{proof}

Corollary \ref{corollary M(K)} is related to a result of Yan \cite{yan}, 
where $\K\subset L(P)$ and $P$ is countably additive.

The representation (\ref{representation}) extends beyond $L(\mu)$. 

\begin{corollary}
\label{corollary general} 
Let $\La\subset\R^\Omega$ be a linear space. A linear functional 
$\phi$ on $\La$ admits a representing measure if and only if there 
exists $\mu\in ba$ such that $\La\subset L(\mu)$ and $\phi$ is 
continuous with respect to the norm topology of $L(\mu)$. If, in 
addition, $\phi$ is positive and $\La$ a vector sublattice of $\R^\Omega$, 
there exists a positive representing measure.
\end{corollary}

\begin{proof}
The direct implcation is obvious. For the converse, let $\mu\in ba$ be as 
in the statement and denote by $\bar\phi$ the continuous, linear extension 
of $\phi$ to $L(\mu)$. If $\La$ is a vector lattice and $\phi$ is positive, 
the inequality $\phi(f)\leq\bar\phi(f^+)$ implies that such extension may
be chosen to be positive and continuous. In either case the claim follows 
from Theorem \ref{theorem L(mu)}.
\end{proof} 

Daniell theorem also follows easily. 

\begin{corollary}
\label{corollary daniell} Let $\La$ be a vector sublattice of $\R^\Omega$ 
containing $1$ and $\phi$ a positive linear functional on $\La$. Then 
$\lim_n\phi(f_n)=0$ for every sequence $\seqn{f}$ in $\La$ which decreases 
to $0$ pointwise if and only if $\phi$ admits a representing measure $m$ 
which is countably additive in restriction to the $\sigma$ algebra 
generated by $\La$.
\end{corollary}

\begin{proof}
Consider the case $\phi\neq 0$, the claim being otherwise trivial. 
Then, by (\ref{riesz}), $\phi(1)>0$ and $\phi$ admits a representing 
probability $m$. Let 
$
\A=\left\{E\subset\Omega:\inf_{\{g\in\La:g\geq \set{E}\}}m(g)=\sup_{\{f\in\La:f\leq\set{E}\}}m(f)\right\}
$
and consider a decreasing sequence $\seqn{E}$ in $\A$ with $\bigcap_nE_n=\varnothing $. 
For each $\eta>0$ there are sequences $\seqn{f}$ and $\seqn{g}$ in $\La_+$ with 
$g_n\geq\set{E_n}\geq f_n$ and $m(f_n)\geq m(g_n) -\eta 2^{-n}$. Let 
$h_n=\inf_{\{k\leq n\}}f_{k}$. $m(h_1)\geq m(g_1) -\eta 2^{-1}$; 
if $m(h_{n-1})\geq m(g_{n-1}) -\eta\sum_{k=1}^{n-1}2^{-k}$ for some $n$ then, 
$h_{n-1}+f_n=h_n+(h_{n-1}\vee f_n)\leq h_n+g_{n-1}$ implies
\begin{eqnarray*}
m(h_n)\geq m(f_n)+m(h_{n-1})-m(g_{n-1})
            \geq m(f_n)-\eta\sum_{k=1}^{n-1}2^{-k}
           \geq m(g_n) -\eta\sum_{k=1}^{n}2^{-k}
\end{eqnarray*}
Thus the sequence $\seqn{f}$ may be chosen to be decreasing to $0$ and such that 
$m(f_n)\geq m(g_n) -\eta $ for each $n$. Then, $0=\lim_nm(f_n)\geq\lim_nm(E_n) -\eta$. 
It is well known that $\A$ is an algebra and that $\La\cap\B\subset\B(\A)$, see e.g. 
\cite[p. 774]{bochner}. Thus, $\restr{m}{\A}$ admits a countably additive extension to 
$\sigma\A$ and this, in turn, an extension $\mu $ to $2^{\Omega }$. Since $\mu$ and 
$m$ coincide on $\A$, $\mu $ is another representing measure for $\phi$. 
The converse is a straightforward implication of monotone convergence.
\end{proof}

\end{document}